\documentclass{article}

\title{\LARGE \textbf{Nonhamiltonian Graphs with Given Toughness }}
\author{Zh.G. Nikoghosyan\footnote{G.G. Nicoghossian (up to 1997)}}

\begin{document}

\maketitle

\begin{abstract}
In 1973, Chv\'{a}tal introduced the concept of toughness $\tau$ of a graph and constructed an infinite class of  nonhamiltonian graphs with $\tau=\frac{3}{2}$. Later Thomassen found nonhamiltonian graphs with $\tau>\frac{3}{2}$, and Enomoto et al. constructed nonhamiltonian graphs with $\tau=2-\epsilon$ for each positive $\epsilon$. The last result in this direction is due to Bauer, Broersma and Veldman, which states that for each positive $\epsilon$, there exists a nonhamiltonian graph with $\tau\ge \frac{9}{4}-\epsilon$. In this paper we prove that for each rational number $t$ with $0<t<\frac{9}{4}$, there exists a nonhamiltonian graph with $\tau=t$.\\


\end{abstract}

\section{Introduction}

Only finite undirected graphs without loops or multiple edges are considered.  The set of vertices of a graph $G$ is denoted by $V(G)$ and the set of edges by $E(G)$.  The order and the independence number of  $G$ is denoted by $n$ and $\alpha$, respectively. For $S$ a subset of $V(G)$, we denote by $G\backslash S$ the maximum subgraph of $G$ with vertex set $V(G)\backslash S$. The neighborhood of a vertex $x\in V(G)$ is denoted by $N(x)$.  A graph $G$ is hamiltonian if $G$ contains a Hamilton cycle, i.e. a cycle of length $n$.  A good reference for any undefined terms is \cite{[5]}. 

The concept of toughness of a graph was introduced in 1973 by Chv\'{a}tal \cite{[6]}. Let $\omega(G)$ denote the number of components of a graph $G$. A graph $G$ is $t$-tough if $|S|\ge t\omega(G\backslash S)$ for every subset $S$ of the vertex set $V(G)$ with $\omega(G\backslash S)>1$. The toughness of $G$, denoted $\tau(G)$, is the maximum value of $t$ for which $G$ is $t$-tough (taking $\tau(K_n)=\infty$ for all $n\ge 1$). By the definition, toughness $\tau$ is a rational number. Since then significant progress has been made toward understanding the relationship between the toughness of a graph and its cycle structure. Much of the research on this subject have been inspired by the following conjecture due to Chv\'{a}tal \cite{[6]}.\\

\noindent\textbf{Conjecture 1}. There exists a finite constant $t_0$ such that every $t_0$-tough graph is hamiltonian.\\

In \cite{[6]}, Chv\'{a}tal constructed an infinite family of nonhamiltonian graphs with $\tau=\frac{3}{2}$, and then Thomassen $[\cite{[4]}, p. 132]$  found nonhamiltonian graphs with $\tau>\frac{3}{2}$. Later Enomoto et al. \cite{[7]} have found nonhamiltonian graphs with $\tau=2-\epsilon$ for each positive $\epsilon$. The last result in this direction is due to Bauer, Broersma and Veldman \cite{[2]} inspired by special constructions introduced in \cite{[1]} and \cite{[3]}.\\

\noindent\textbf{Theorem A}. For each positive $\epsilon>0$, there exists a nonhamiltonian graph with $\frac{9}{4}-\epsilon<\tau<\frac{9}{4}$.\\

In view of Theorem A, the following problem seems quite reasonable.\\

\noindent\textbf{Problem}. Is there a nonhamiltonian graph $G$ with $\tau(G)=t$ for a given rational number $t$ with $0<t< \frac{9}{4}$?\\

In this paper we prove the following.\\

\noindent\textbf{Theorem 1}. For each rational number $t$ with $0<t<\frac{9}{4}$, there exists a nonhamiltonian graph $G$ with $\tau(G)=t$.\\

\section{Preliminaries}

To prove Theorem 1, we need the following graph constructions. \\

\textbf{Definition 1}. Let $L^{(1)}$ be a graph obtained from $C_8=w_1w_2...w_8w_1$ by adding the edges $w_2w_4, w_4w_6, w_6w_8$ and $w_2w_8$. Put $x=w_1$ and $y=w_5$. This is the well-known building block $L$ used to obtain $(\frac{9}{4}-\epsilon)$-tough nonhamiltonian graphs (see \cite{[2]}, Figure 1).   \\

In this paper we will use a number of additional modified building blocks.\\

\textbf{Definition 2}. Let $L^{(2)}$  be the graph obtained from $L^{(1)}$ by deleting the edges $w_1w_2$, $w_2w_8$ and identifying $w_2$ with $w_8$. \\

\textbf{Definition 3}. Let $L^{(3)}$ be the graph obtained from $L^{(1)}$ by adding  a new vertex $w_9$ and the edges $w_4w_9$, $w_6w_9$.\\

\textbf{Definition 4}. Let $L^{(4)}$ be the graph obtained from the triangle $w_1w_2w_3w_1$ by adding the vertices $w_4, w_5$ and the edges $w_1w_4$, $w_3w_5$. Put $x=w_4$ and $y=w_5$.\\

\textbf{Definition 5}. For each $L\in \{L^{(1)}, L^{(2)}\}$, define the graph  
$G(L,x,y,l,m)$ $(l,m\in N)$ as follows. Take $m$ disjoint copies 
$L_1,L_2,...,L_m$ of $L$, with $x_i,y_i$ the vertices 
in $L_i$ corresponding to the vertices $x$ and $y$ in 
$L$ $(i=1,2,...,m)$. Let $F_m$ be the graph obtained 
from $L_1\cup ... \cup L_m$ by adding all possible edges
 between pairs of vertices in $x_1,...,x_m,y_1,...,y_m$. 
Let $T=K_l$ and let $G(L,x,y,l,m)$ be the join $T\vee F_m$ of $T$ and $F_m$. \\

The following can be checked easily.\\

\textbf{Claim 1}. The vertices $x$ and $y$ are not connected by a Hamilton path of $L^{(i)}$ $(i=1,2,3)$.\\

The proof of the following result occurs in \cite{[1]}, which we repeat here for convenience.\\

\textbf{Claim 2}. Let $H$ be a graph and $x,y$ two vertices of $H$ which are not connected by a Hamilton path of $H$. If $m\ge 2l+1$ then $G(H,x,y,l,m)$ is nonhamiltonian.\\

\textbf{Proof}. Suppose $G(H,x,y,l,m)$ contains a 
Hamilton cycle $C$. The intersection of $C$ and 
$F_m$ consists of a collection $\Re$ of at most 
$l$ disjoint paths, together containing all vertices 
in $F_m$. Since $m\ge 2l+1$, there is a subgraph 
$H_{i_0}$ in $F_m$ such that no endvertex of a path 
of $\Re$ lies in $H_{i_0}$. Hence the intersection of 
$C$ and $H_{i_0}$ is a path with endvertices
 $x_{i_0}$ and $y_{i_0}$ that contains all vertices
 of $H_{i_0}$. This contradicts the fact that $H_{i_0}$ is 
a copy of the graph $H$ without a Hamilton path
 between $x$ and $y$. Claim 1 is proved.\\

 \section{Proof of Theorem 1}

Let $t$ be a rational number with $0<t<\frac{9}{4}$ and let $t=\frac{a}{b}$ for some integers $a,b$. \\

\textbf{Case 1}. $0<\frac{a}{b}<1$.

Let $K_{a,b}$  be the complete bipartite graph $G=(V_1,V_2;E)$ with vertex classes $V_1$ and $V_2$ of order $a$ and $b$, respectively. Since $\frac{a}{b} <1$,  we have $\alpha(G)=b>(a+b)/2$ and therefore, $K_{a,b}$ is a nonhamiltonian graph. Clearly, $\tau\le |V_1|/\omega(G\backslash V_1)=a/b$. Choose $S\subset V(G)$ such that $\tau(K_{a,b})=|S|/\omega(G\backslash S)$. Put $S\cap V_i=S_i$ and $|S_i|=s_i$ $(i=1,2)$. If $V_i\backslash S\not=\emptyset$ $(i=1,2)$ then clearly $\omega(G\backslash S)=1$, which is impossible by the definition. Hence $V_i\backslash S=\emptyset$ for some $i\in \{1,2\}$, i.e. $V_i\subseteq S$.\\

\textbf{Case 1.1}. $i=2$.

Since $\tau=(s_1+b)/(a-s_1)\ge b/a$, we have $s_1=0$, i.e. $S=V_2$ and $\tau=b/a$, contradicting that fact that $\tau\le a/b$. \\

\textbf{Case 1.2}. $i=1$.

Since $\tau=(s_2+a)/(b-s_2)\ge a/b$, we have $s_2=0$, implying that $S=V_1$ and $\tau=a/b$.\\

\textbf{Case 2}. $\frac{a}{b}=1$.

Let $G$ be a graph obtained from $C_6=x_1x_2...x_6x_1$ by adding a new vertex $x_7$ and the edges $x_1x_7,x_4x_7, x_2x_6$. Clearly, $G$ is not hamiltonian and $\tau(G)=1$.\\

\textbf{Case 3}. $1<\frac{a}{b}<\frac{3}{2}$.

\textbf{Case 3.1}. $\frac{a}{b}<\frac{3}{2}-\frac{1}{b}$.

Let $V_1,V_2,V_3$ be pairwise disjoint sets of vertices with 
$$
V_1=\{x_1,x_2,...,x_{a-b+1}\}, \ V_2=\{y_1,y_2,...,y_b\}, \ V_3=\{z_1,z_2,...,z_b\}.
$$
Join each $x_i$ to all the other vertices and each $z_i$ to every other $z_j$ as well as to the vertex $y_i$ with the same subscript $i$. Call the resulting graph $H$. To determine the toughness of $H$, choose $W\subset V(H)$ such that $\tau(H)=|W|/\omega(H\backslash W)$. Put $m=|W\cap V_3|$. Clearly, $W$ is a minimal set whose removal from $H$ results in a graph with $\omega(H\backslash W)$ components. As $W$ is a cutset, we have $V_1\subset W$ and $m\ge1$. From the minimality of $W$ we easily conclude that $V_2\cap W=\emptyset$ and $m\le b-1$. Then we have $|W|=m+a-b+1$ and $\omega(H\backslash W)=m+1$. Hence
$$
\tau(H)=\frac{|W|}{\omega(H\backslash W)}=\min_{1\le m\le b-1}\frac{m+a-b+1}{m+1}=\frac{a}{b}.
$$
To see that $H$ is nonhamiltonian, let us assume the contrary, i.e. let $C$ be a Hamilton cycle in $H$. Denote by $F$ the set of edges of $C$ having at least one endvertex in $V_2$. Since $V_2$ is independent, we have $|F|=2|V_2|$. On the other hand, there are at most $2|V_1|$ edges in $F$ having one endvertex in $V_1$ and at most $|V_3|$ edges in $F$ having one endvertex in $V_3$. Thus
$$
2b=2|V_2|=|F|\le2|V_1|+|V_3|=2(a-b+1)+b=2a-b+2.
$$
But this is equivalent to $a/b\ge 3/2-1/b$, contradicting the hypothesis.\\

\textbf{Case 3.2}. $\frac{a}{b}\ge\frac{3}{2}-\frac{1}{b}$.

By choosing $q\in N$ sufficiently large with 
$$
\frac{a}{b}=\frac{aq}{bq}<\frac{3}{2}-\frac{1}{bq},
$$
we can argue as in Case 3.1.\\

\textbf{Case 4}. $\frac{a}{b}=\frac{3}{2}$.

An example of a nonhamiltonian graph with $\tau=3/2$ is obtained when in the Petersen graph,
each vertex is replaced by a triangle.\\

\textbf{Case 5}. $\frac{3}{2}<\frac{a}{b}<\frac{7}{4}$.\\

\textbf{Claim 3}. For $l\ge 2$ and $m\ge1$,
$$
\tau(G(L^{(2)},x,y,l,m))=\frac{l+3m}{1+2m}.
$$ 
\textbf{Proof}. Let $G=G(L^{(2)},x,y,l,m)$ for some $l\ge2$ and $m\ge1$. Choose $S\subseteq V(G)$ such that $\omega(G\backslash S)>1$ and $\tau(G)=|S|/\omega(G\backslash S)$. Obviously, $V(T)\subseteq S$. Define $S_i=S\cap V(L_i)$, $s_i=|S_i|$, and let $\omega_i$ be the number of components of $L_i\backslash S_i$ that contain neither $x_i$ nor $y_i$ $(i=1,...,m)$. Then
$$
\tau(G)=\frac{l+\sum_{i=1}^ms_i}{c+\sum_{i=1}^m\omega_i}\ge \frac{l+\sum_{i=1}^ms_i}{1+\sum_{i=1}^m\omega_i},
$$
where $c=0$ if $x_i,y_i\in S$ for all $i\in \{1,...,m\}$ and $c=1$ otherwise. It is easy to see that
$$
\omega_i\le2, \  \  s_i\ge \frac{3}{2}\omega_i  \  \  (i=1,...,m).
$$
Then
$$
\tau\ge \frac{l+\frac{3}{2}\sum_{i=1}^m\omega_i}{1+\sum_{i=1}^m\omega_i}=\frac{l-\frac{3}{2}}{1+\sum_{i=1}^m\omega_i}+\frac{3}{2}
$$
$$
\ge\frac{l-\frac{3}{2}}{1+2m}+\frac{3}{2}=\frac{l+3m}{1+2m}.
$$

Set $U=V(T)\cup U_1\cup ...\cup U_m$, where $U_i$ is the set of vertices of $L_i$ having degree at least 4 in $L_i$ $(i=1,...,m)$. The proof of Claim 3 is completed by observing that 
$$
\tau(G)\le\frac{|U|}{\omega(G\backslash U)}=\frac{l+3m}{2m+1}.           \quad     \quad     \rule{7pt}{6pt} 
$$
 
\textbf{Case 5.1}. $b=2k+1$ for some integer $k$.

Consider the graph $G(L^{(2)},x,y,a-\frac{3}{2}(b-1),\frac{b-1}{2})$. \\

\textbf{Case 5.1.1}. $\frac{a}{b}\le\frac{7}{4}-\frac{9}{4b}$.

By the hypothesis,
$$
m=\frac{b-1}{2}\ge2(a-\frac{3}{2}(b-1))+1=2l+1,
$$
implying by Claim 2 that $G$ is not hamiltonian. Clearly $b\ge3$, implying that $m=(b-1)/2\ge1$. \\

\textbf{Case 5.1.1.1}. $\frac{a}{b}\ge\frac{3}{2}+\frac{1}{2b}$.

By the hypothesis, $l=a-\frac{3}{2}(b-1)\ge2$. By Claim 3, $\tau(G)=\frac{a}{b}$.\\

\textbf{Case 5.1.1.2}. $\frac{a}{b}<\frac{3}{2}+\frac{1}{2b}$.

By choosing a sufficiently large integer $q$ with 
$$
\frac{a}{b}=\frac{aq}{bq}\ge\frac{3}{2}+\frac{1}{2bq},
$$
we can argue as in Case 5.1.1.1.\\

\textbf{Case 5.1.2}. $\frac{a}{b}>\frac{7}{4}-\frac{9}{4b}$.

By choosing a sufficiently large integer $q$ with
$$
\frac{a}{b}=\frac{aq}{bq}\le\frac{7}{4}-\frac{9}{4bq},
$$
we can argue as in Case 5.1.1.\\

\textbf{Case 5.2}. $b=2k$ for some integer $k$.

Consider the graph $G^\prime$ obtained from $G(L^{(2)},x,y,l,m)$ by replacing $L_m$ with $L^{(3)}$. \\

\textbf{Claim 4}. For $l\ge 2$ and $m\ge1$,
$$
\tau(G^\prime)=\frac{l+3m+1}{2(m+1)}.
$$ 

\textbf{Proof}. Choose $S\subseteq V(G^\prime)$ such that $\omega(G^\prime\backslash S)>1$ and $\tau(G^\prime)=|S|/\omega(G^\prime\backslash S)$. Obviously, $V(T)\subseteq S$. Define $S_i=S\cap V(L_i)$, $s_i=|S_i|$, and let $\omega_i$ be the number of components of $L_i\backslash S_i$ that contain neither $x_i$ nor $y_i$ $(i=1,...,m)$. Since $s_i\ge\frac{3}{2}\omega_i$ $(i=1,...,m-1)$ and $s_m\ge\frac{4}{3}\omega_m$, we have
$$
\tau(G^\prime)\ge \frac{l+\sum_{i=1}^ms_i}{c+\sum_{i=1}^m\omega_i}\ge\frac{l+\frac{3}{2}\sum_{i=1}^{m-1}\omega_i+\frac{4}{3}\omega_m}{1+\sum_{i=1}^m\omega_i}=\frac{l-\frac{1}{6}\omega_m}{1+\sum_{i=1}^m\omega_i}+\frac{3}{2},
$$
where $c=0$ if $x_i,y_i\in S$ for all $i\in \{1,...,m\}$ and $c=1$ otherwise. Observing also that $\omega_i\le2$ $(i=1,...,m-1)$ and $\omega_m\le3$, we obtain
$$
(l-2)\sum_{i=1}^m\omega_i+\frac{1}{3}(m+1)\omega_m\le(l-2)(2m+1)+(m+1)\le2l(m+1).
$$ 
But this is equivalent to 
$$
\frac{l-\frac{1}{6}\omega_m}{1+\sum_{i=1}^m\omega_i}+\frac{3}{2}\ge\frac{l-2}{2(m+1)}+\frac{3}{2},
$$
implying that 
$$
\tau(G^\prime)\ge\frac{l-2}{2(m+1)}+\frac{3}{2}=\frac{l+3m+1}{2(m+1)}.
$$

Set $U=V(T)\cup U_1\cup ...\cup U_m$, where $U_i$ is the set of vertices of $L_i$ having degree at least 4 in $L_i$ $(i=1,...,m)$. The proof of Claim 4 is completed by observing that 
$$
\tau(G^\prime)\le\frac{|U|}{\omega(G\backslash U)}=\frac{l+3m+1}{2(m+1)}.           \quad     \quad     \rule{7pt}{6pt} 
$$
Consider the graph $G^\prime$ with $m=\frac{b}{2}-1$ and $l=a-\frac{3}{2}b+2$. Clearly $m=\frac{b}{2}-1\ge1$ and $l=a-\frac{3}{2}b+2\ge2$. By Claim 4, $\tau(G^\prime)=\frac{a}{b}$.   \\

\textbf{Case 5.2.1}. $\frac{a}{b}\le\frac{7}{4}-\frac{3}{b}$.

By the hypothesis, $m\ge2l+1$, and by Claim 2, $G^\prime$ is not hamiltonian. \\

\textbf{Case 5.2.2}. $\frac{a}{b}>\frac{7}{4}-\frac{3}{b}$.

By choosing a sufficiently large $q$ with 
$$
\frac{a}{b}=\frac{aq}{bq}\le\frac{7}{4}-\frac{3}{b},
$$
we can argue as in Case 5.2.1.\\

\textbf{Case 6}. $\frac{7}{4}-\epsilon<\frac{a}{b}\le2$.

Let $m=m_1+m_2\ge2l+1$ and let $G^{\prime\prime}$ be the graph obtained from $G(L^{(1)},x,y,l,m)$ by replacing $L_i$ with $L^{(2)}$ $(i=m_1+1,m_1+2,...,m)$. By Claim 2, $G^{\prime\prime}$ is not hamiltonian.   \\

\textbf{Claim 5}. For $l\ge 2$, $m\ge1$ and $m_2\ge l-2$,
$$
\tau(G^{\prime\prime})=\frac{l+3m_2}{2m_2+1}.
$$ 

\textbf{Proof}. Choose $S\subseteq V(G^{\prime\prime})$ such that 
$$
\omega(G^{\prime\prime}\backslash S)>1,   \  \  \tau(G^{\prime\prime})=|S|/\omega(G^{\prime\prime}\backslash S).
$$
Obviously, $V(T)\subseteq S$. Define $S_i=S\cap V(L_i)$, $s_i=|S_i|$, and let $\omega_i$ be the number of components of $L_i\backslash S_i$ that contain neither $x_i$ nor $y_i$ $(i=1,...,m)$. Since $s_i\ge2\omega_i$ $(i=1,...,m_1)$ and $s_i\ge\frac{3}{2}\omega_i$ $(i=m_1+1,...,m)$, we have
$$
\tau(G^{\prime\prime})\ge \frac{l+\sum_{i=1}^{m_1}s_i+\sum_{i=m_1+1}^ms_i}{c+\sum_{i=1}^m\omega_i}\ge\frac{l+2\sum_{i=1}^{m_1}\omega_i+\frac{3}{2}\sum_{i=m_1+1}^m\omega_i}{1+\sum_{i=1}^m\omega_i}
$$
$$
\frac{l+\frac{1}{2}\sum_{i=1}^{m_1}\omega_i-\frac{3}{2}+\frac{3}{2}(1+\sum_{i=1}^m\omega_i)}{1+\sum_{i=1}^m\omega_i}=\frac{2l+\sum_{i=1}^{m_1}\omega_i-3}{2(1+\sum_{i=1}^m\omega_i)}+\frac{3}{2},
$$
where $c=0$ if $x_i,y_i\in S$ for all $i\in \{1,...,m\}$ and $c=1$ otherwise. Observing also that $\omega_i\le2$ $(i=1,...,m)$, we obtain 
$$
(2l-3)\sum_{i=m_1+1}^m\omega_i-(2m_2-2l+4)\sum_{i=1}^{m_1}\omega_i\le4lm_2-6m_2.
$$
But this is equivalent to
$$
\frac{2l+\sum_{i=1}^{m_1}\omega_i-3}{2(1+\sum_{i=1}^m\omega_i)}+\frac{3}{2}\ge\frac{2l-3}{2(2m_2+1)}+\frac{3}{2},
$$
implying that
$$
\tau(G^{\prime\prime})\ge \frac{2l-3}{2(2m_2+1)}+\frac{3}{2}=\frac{l+3m_2}{2m_2+1}.
$$
Set $U=V(T)\cup U_1\cup ...\cup U_m$, where $U_i$ is the set of vertices of $L_i$ having degree at least 4 in $L_i$ $(i=1,...,m)$. The proof of Claim 5 is completed by observing that 
$$
\tau(G^{\prime\prime})\le\frac{|U|}{\omega(G\backslash U)}=\frac{l+3m_2}{2m_2+1}.           \quad     \quad     \rule{7pt}{6pt} 
$$

\textbf{Case 6.1}. $b=2k+1$ for some integer $k$.

Consider the graph $G^{\prime\prime}$ with $m_2=\frac{b-1}{2}$ and $l=a-\frac{3}{2}(b-1)$.\\

\textbf{Case 6.1.1}. $\frac{a}{b}\ge\frac{3}{2}+\frac{1}{2b}$.

Since $\frac{a}{b}\le2$, we have 
$$
m_2=\frac{b-1}{2}\ge a-\frac{3}{2}(b-1)-2=l-2.
$$
Next, since $\frac{a}{b}\ge\frac{3}{2}+\frac{1}{2b}$, we have $l=a-\frac{3}{2}(b-1)\ge2$. By Claim 5, $\tau(G^{\prime\prime})=\frac{a}{b}$.\\

\textbf{Case 6.1.2}. $\frac{a}{b}<\frac{3}{2}+\frac{1}{2b}$.

By choosing a sufficiently large integer $q$ with 
$$
\frac{a}{b}=\frac{aq}{bq}\ge\frac{3}{2}+\frac{1}{2bq},
$$
we can argue as in Case 6.1.1.\\

\textbf{Case 6.2}. $b=2k$ for some integer $k$.

Consider the graph $G^{\prime\prime\prime}$ obtained from $G^{\prime\prime}$ by replacing $L_m$ with $L^{(3)}$. \\

\textbf{Claim 6}. For $l\ge 2$, $m\ge1$ and $m_2\ge l-2$,
$$
\tau(G^{\prime\prime\prime})=\frac{l+3m_2+1}{2(m_2+1)}.
$$ 

\textbf{Proof}. Choose $S\subseteq V(G^{\prime\prime\prime})$ such that 
$$
\omega(G^{\prime\prime\prime}\backslash S)>1, \  \  \tau(G^{\prime\prime\prime})=|S|/\omega(G^{\prime\prime\prime}\backslash S)
$$
Obviously, $V(T)\subseteq S$. Define $S_i=S\cap V(L_i)$, $s_i=|S_i|$, and let $\omega_i$ be the number of components of $L_i\backslash S_i$ that contain neither $x_i$ nor $y_i$ $(i=1,...,m)$. Since $s_i\ge2\omega_i$ $(i=1,...,m_1)$, $s_i\ge\frac{3}{2}\omega_i$ $(i=m_1+1,...,m-1)$ and $s_m\ge\frac{4}{3}\omega_m$, we have

$$
\tau(G^{\prime\prime\prime})\ge \frac{l+\sum_{i=1}^{m_1}s_i+\sum_{i=m_1+1}^{m-1}s_i+s_m}{c+\sum_{i=1}^m\omega_i}
$$
$$
\ge\frac{l+2\sum_{i=1}^{m_1}\omega_i+\frac{3}{2}\sum_{i=m_1+1}^{m-1}\omega_i+\frac{4}{3}\omega_m}{1+\sum_{i=1}^m\omega_i}
$$
$$
=\frac{l+\frac{1}{2}\sum_{i=1}^{m_1}\omega_i-\frac{1}{6}\omega_m+(\frac{3}{2}\sum_{i=1}^{m_1}\omega_i+\frac{3}{2}\sum_{i=m_1+1}^{m})}{1+\sum_{i=1}^m\omega_i}
$$
$$
=\frac{l+\frac{1}{2}\sum_{i=1}^{m_1}\omega_i-\frac{1}{6}\omega_m}{1+\sum_{i=1}^m\omega_i}+\frac{3}{2},
$$
where $c=0$ if $x_i,y_i\in S$ for all $i\in \{1,...,m\}$ and $c=1$ otherwise. Observing also that $\omega_i\le2$ $(i=1,...,m-1)$ and $\omega_m\le3$, we obtain 
$$
(l-2)\sum_{i=m_1+1}^{m}\omega_i+\frac{1}{3}(m_2+1)\omega_m-(m_2-l+3)\sum_{i=1}^{m_1}\omega_i\le l+2lm_2+2.
$$
But this is equivalent to
$$
\frac{l+\frac{1}{2}\sum_{i=1}^{m_1}\omega_i-\frac{1}{6}\omega_m}{1+\sum_{i=1}^m\omega_i}+\frac{3}{2}\ge \frac{l-2}{2(m_2+1)}+\frac{3}{2},
$$
implying that
$$
\tau(G^{\prime\prime\prime})\ge \frac{l-2}{2(m_2+1)}+\frac{3}{2}=\frac{l+3m_2+1}{2(m_2+1)}.
$$
Set $U=V(T)\cup U_1\cup ...\cup U_m$, where $U_i$ is the set of vertices of $L_i$ having degree at least 4 in $L_i$ $(i=1,...,m)$. The proof of Claim 6 is completed by observing that 
$$
\tau(G^{\prime\prime\prime})\le\frac{|U|}{\omega(G\backslash U)}=\frac{l+3m_2+1}{2(m_2+1)}.           \quad     \quad     \rule{7pt}{6pt} 
$$
Consider the graph $G^{\prime\prime\prime}$ with $m_2=\frac{b}{2}-1$ and $l=a-\frac{3}{2}b+2$.\\

\textbf{Case 6.2.1}. $\frac{a}{b}\le2-\frac{1}{b}$.

By the hypothesis, $m_2=\frac{b}{2}-1\ge(a-\frac{3}{2}b+2)-2=l-2$. Next, since $\frac{a}{b}>\frac{7}{4}-\epsilon>\frac{3}{2}$, we have $l=\frac{3}{2}b+2\ge2$. By Claim 6, $\tau(G^{\prime\prime\prime})=\frac{a}{b}$.\\

\textbf{Case 6.2.2}. $\frac{a}{b}>2-\frac{1}{b}$.

By choosing a sufficiently large integer $q$ with $\frac{a}{b}=\frac{aq}{bq}\le2-\frac{1}{bq}$, we can argue as in Case 6.2.1.\\

\textbf{Case 7}. $2<\frac{a}{b}<\frac{9}{4}$.

\textbf{Case 7.1}. $b=2k+1$ for some integer $k$.

\textbf{Case 7.1.1}. $\frac{a}{b}\le\frac{9}{4}-\frac{11}{4b}$.

Take the graph $G(L^{(1)},x,y, a-2b+2,\frac{b-1}{2})$. Since $\frac{a}{b}>2$, we have $l=a-2b+2\ge2$. Next, the hypothesis $\frac{a}{b}\le\frac{9}{4}-\frac{11}{4b}$ is equivalent to 
$$
m=\frac{b-1}{2}\ge2(a-2b+2)+1=2l+1.
$$
By Claim 1,  $G(L^{(1)},x,y, a-2b+2,\frac{b-1}{2})$ is not hamiltonian. The toughness $\tau(G(L^{(1)},x,y, a-2b+2,\frac{b-1}{2}))$ can be determined exactly as in proof of Theorem A \cite{[2]},
$$
\tau(G(L^{(1)},x,y, a-2b+2,\frac{b-1}{2}))\ge\frac{l+4m}{2m+1}=\frac{a}{b}. 
$$

\textbf{Case 7.1.2}. $\frac{a}{b}>\frac{9}{4}-\frac{11}{4b}$.

By choosing a sufficiently large integer $q$ with 
$$
\frac{aq}{bq}=\frac{a}{b}\le\frac{9}{4}-\frac{11}{4bq},
$$
we can argue as in Case 7.1.1.\\

\textbf{Case 7.2}. $b=2k$ for some positive integer $k$.

Take the graph $G^{\prime\prime\prime\prime}$ obtained from $G(L^{(1)},x,y, a-2b+2,\frac{b}{2})$ by replacing $L_m$ with $L^{(4)}$. Since $\frac{a}{b}>2$, we have $l=a-2b+2>2$. We have also $m=\frac{b}{2}>1$, since $b\ge3$.\\

\textbf{Claim 7}. For $l\ge 2$ and  $m\ge1$,
$$
\tau(G^{\prime\prime\prime\prime})=\frac{l+4m-2}{2m}.
$$ 

\textbf{Proof}. Choose $S\subseteq V(G^{\prime\prime\prime})$ such that 
$$
\omega(G^{\prime\prime\prime}\backslash S)>1, \  \  \tau(G^{\prime\prime\prime})=|S|/\omega(G^{\prime\prime\prime}\backslash S)
$$
Obviously, $V(T)\subseteq S$. Define $S_i=S\cap V(L_i)$, $s_i=|S_i|$, and let $\omega_i$ be the number of components of $L_i\backslash S_i$ that contain neither $x_i$ nor $y_i$ $(i=1,...,m)$. Since $s_i\ge2\omega_i$ $(i=1,...,m)$, $\omega_i\le2$ $(i=1,...,m-1)$ and $\omega_m\le1$, we have
$$
\tau(G^{\prime\prime\prime\prime})=\frac{l+\sum_{i=1}^ms_i}{c+\sum_{i=1}^m\omega_i}\ge\frac{l+2\sum_{i=1}^m\omega_i}{1+\sum_{i=1}^m\omega_i}
$$
$$
=\frac{l-2}{1+\sum_{i=1}^m\omega_i}+2\ge \frac{l-2}{2m}+2=\frac{l+4m-2}{2m},
$$
where $c=0$ if $x_i,y_i\in S$ for all $i\in \{1,...,m\}$ and $c=1$ otherwise. Set $U=V(T)\cup U_1\cup ...\cup U_m$, where $U_i$ is the set of vertices of $L_i$ having degree at least 4 in $L_i$ $(i=1,...,m)$. The proof of Claim 7 is completed by observing that 
$$
\tau(G^{\prime\prime\prime\prime})\le\frac{|U|}{\omega(G\backslash U)}=\frac{l+4m-2}{2m}.           \quad     \quad     \rule{7pt}{6pt} 
$$

\textbf{Case 7.2.1}. $\frac{a}{b}\le\frac{9}{4}-\frac{3}{b}$.

By the hypothesis,
$$
m-1=\frac{b}{2}-1\ge2(a-2b+2)+1=2l+1.
$$
By Claim 2, $G^{\prime\prime\prime\prime}$ is not hamiltonian. By Claim 7, $\tau(G^{\prime\prime\prime\prime})=\frac{a}{b}$.\\

\textbf{Case 7.2.2}. $\frac{a}{b}>\frac{9}{4}-\frac{3}{b}$.

By choosing a sufficiently large integer $q$ with
$$
\frac{aq}{bq}=\frac{a}{b}\le\frac{9}{4}-\frac{3}{3bq},
$$
we can argue as in Case 7.2.1. Theorem 1 is proved.     \quad   \quad         \rule{7pt}{6pt} \\

\noindent Institute for Informatics and Automation Problems\\ National Academy of Sciences\\
P. Sevak 1, Yerevan 0014, Armenia\\ E-mail: zhora@ipia.sci.am

\end{document}